\documentclass[12pt]{amsart}
\usepackage{amsthm}
\usepackage{amssymb}
\usepackage{amsmath}
\usepackage{graphicx,multicol} 
\usepackage{color}
\theoremstyle{plain}
\newtheorem{proposition}{Proposition}
\newtheorem{theorem}[proposition]{Theorem}
\newtheorem{lemma}[proposition]{Lemma}

\theoremstyle{definition}

\theoremstyle{definition}

\newtheorem{remark}[proposition]{Remark}

\numberwithin{equation}{section}

\numberwithin{proposition}{section}

\gdef\myletter{}

\let\savetheequation\theequation
\def\theequation{\savetheequation\myletter}

\usepackage{color}

\newcommand{\CC}{{\mathbb C}}

\newcommand{\RR}{{\mathbb R}}
\newcommand{\ZZ}{{\mathbb Z}}

\renewcommand{\date}{\today}

\begin{document}


\title[Pac-Man]{\bf Monge-Amp\`ere of Pac-Man}

\author{Norm Levenberg* and Sione Ma'u}{\thanks{*Supported by Simons Foundation grant No. 354549}}

\maketitle
\begin{abstract} We show that the Monge-Amp\`ere density of the extremal function $V_P$ for a non-convex Pac-Man set $P\subset \RR^2$ tends to a finite limit as we approach the vertex $p$ of $P$ linearly but with a value that may vary with the line. On the other hand, along a tangential approach to $p$ the Monge-Amp\`ere density becomes unbounded. This partially  mimics the behavior of the Monge-Amp\`ere density of the union of two quarter disks set $S$ of Sigurdsson and Snaebjarnarson \cite{SS}. We also recover their formula for $V_S$ by elementary methods.
\end{abstract}

\section{Introduction} Given a compact set $K\subset \CC^d$, let
$$V_K(z) := \sup\{u(z):u\in L(\CC^d), \ u\leq 0 \ \hbox{on} \ K\}$$
where $L(\CC^d):= \{u\in PSH(\CC^d): u(z)- \log|z| =0(1), \ |z| \to \infty \}$ and $PSH(\CC^d)$ denotes the cone of plurisubharmonic functions on $\CC^d$. The Siciak-Zaharjuta extremal function is defined as
$$V_K^*(z):=\limsup_{\zeta \to z}V_K(\zeta).$$ Moreover, one has 
$$V_K(z)=\max[0,\sup  \{  \frac{1}{deg (p)}\log |p(z)|: ||p||_K:=\max_{\zeta \in K}|p(\zeta)| \leq 1\}]$$ where the supremum is taken over all non-constant holomorphic polynomials $p$. If $K$ is not pluripolar; i.e., for any $u\in PSH(\CC^d)$ with $u|_K=-\infty$ we have $u\equiv -\infty$, then $V_K^*\in L(\CC^d)$ and the Monge-Amp\`ere measure 
$$\mu_K:=(dd^c V^*_K)^d,$$
which is supported in $K$, plays the role, if $d>1$, of the potential-theoretic equilibrium measure from the case when $d=1$.

Explicit calculation of $V_K^*$ and $\mu_K$ is difficult in general; however qualitative properties of $\mu_K$ are known in certain special cases. Indeed, if $K\subset \RR^d\subset \CC^d$ is a convex body in $\RR^d$; i.e., a compact, convex set with non-empty interior, much is known. Going back to the work of Bedford and Taylor in \cite{BT}, for such $K\subset \RR^d 
$ with smooth boundary, $\mu_K$ is absolutely continuous with respect to Lebesgue measure on $\RR^d$ with density $\rho_K(x)$ which is comparable to $ [dist (x,\partial K)]^{-1/2}$ (Theorem 1.2 of \cite{BT}). In particular, this density blows up as $x$ approaches the boundary of $K$. More generally, for $K\subset \RR^d 
$ compact, on compact subsets of the interior of $K$ the Monge-Amp\`ere measure $\mu_K$ is equivalent to Lebesgue measure on $\RR^d$ (see Theorem 1.1 of \cite{BT} for a precise statement). 

In order simply for $\rho_K(x)\to \infty$ as $x\to \partial K$ the hypothesis of smooth boundary in Theorem 1.2 of \cite{BT} is not necessary (see section 3) but it was unknown to what extent the convexity hypothesis is needed for this conclusion. A few years ago Robert Berman asked what happens to the Monge-Amp\`ere density $\rho_K$ of a non-convex body $K\subset \RR^d$ near a non-convex boundary point such as a fattened $L-$shape in $\RR^2$ at its inner bend or a ``Pac-Man'' (a disk minus a centrally symmetric wedge of opening less than $\pi$) at the vertex of its mouth (Figure 1).

\begin{figure}
\includegraphics[height=35mm]{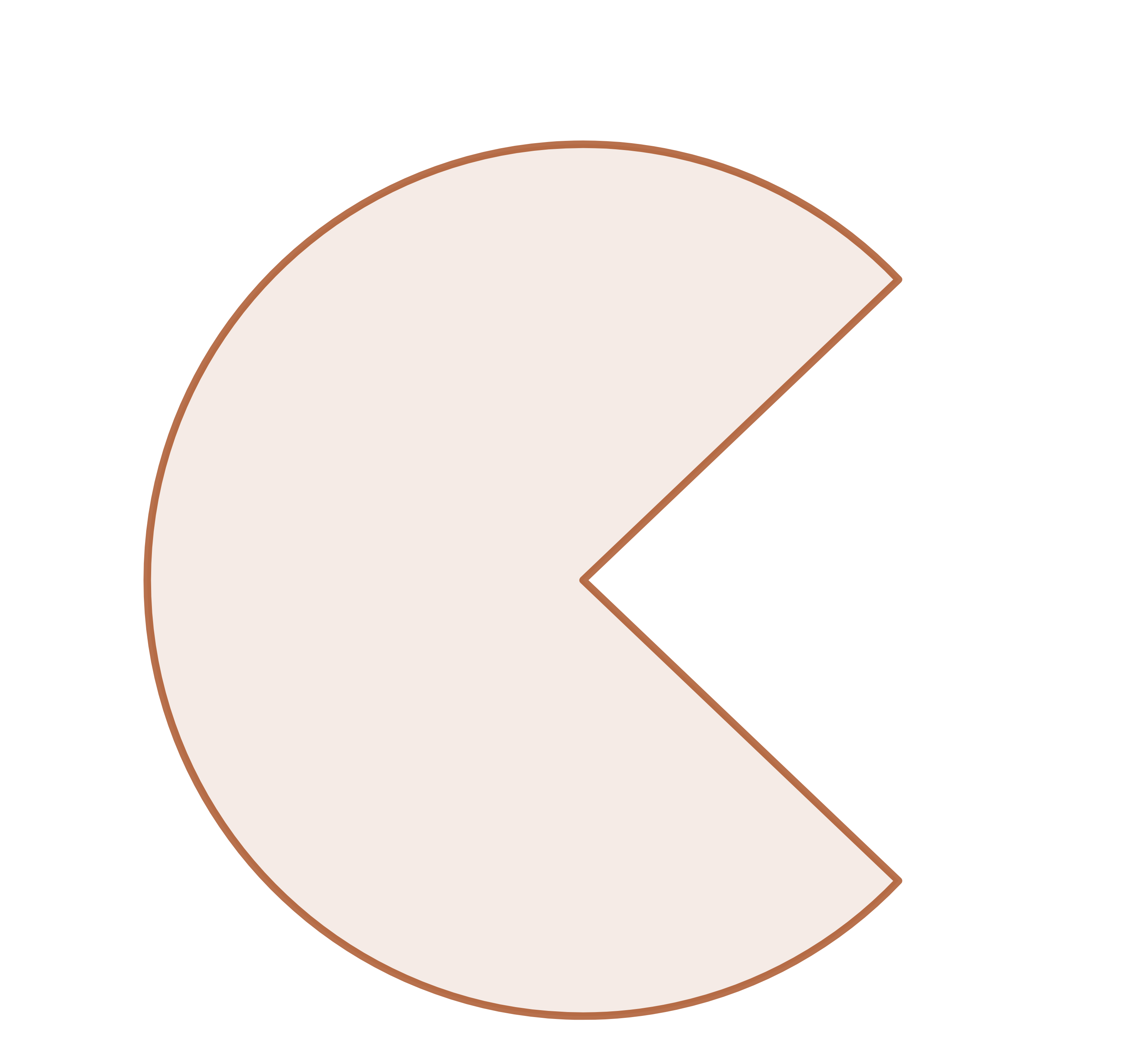}
\caption{A ``Pac-Man'' set}
\end{figure}

 Motivated by an example in a recent preprint of Sigurdsson and Snaebjarnarson \cite{SS}, we show that in these cases, an interesting phenomenon occurs: the Monge-Amp\`ere density of the extremal function tends to a finite limit as we approach the bend or vertex $p$ linearly within $K$ but with a value that varies with the line if the line is sufficiently close to the boundary edge. On the other hand, along a tangential approach to $p$ the Monge-Amp\`ere density becomes unbounded.

In the next section, we recall some results from \cite{BT}, \cite{K} and \cite{SS}. In particular, we provide an elementary proof of the formula for the extremal function of Example 4.4 in \cite{SS}. We give the details of the previously described behavior of the Monge-Amp\`ere density of the extremal function of a Pac-Man set near its vertex in section 3. Finally, in section 4, we indicate the genesis of our method of attack on this problem and discuss directions for future research.

We would like to thank Akil Narayan for pointing out an error in the original statement of our main result, Theorem \ref{pac}.

\section{Background and preliminaries}

A key result which we use is Theorem 5.3.1 of \cite{K}.

\begin{theorem} \label{klimek} Let $F:\CC^d\to \CC^d$ be a proper, polynomial mapping of degree at most $b$ with
$$\liminf_{|z|\to \infty}\frac{||F(z)||}{||z||^a}>0$$
where $b\geq a>0$. Then for any $K\subset \CC^d$ compact,
$$aV_{F^{-1}(K)} \leq V_K(F(z))\leq bV_{F^{-1}(K)}.$$

\end{theorem}

To compare Monge-Amp\`ere measures of extremal functions associated to compact subsets of $\RR^d$, we recall Lemma 2.1 of \cite{BT}.

\begin{lemma} \label{twoone} Let $\Omega$ be a domain in $\CC^d$ and let $u_1,u_2\in PSH(\Omega)\cap L^{\infty}_{loc}(\Omega)$. Suppose $S\subset \Omega\cap \RR^d$ is a closed set containing the supports of the Monge-Amp\`ere measures $(dd^cu_1)^d, (dd^cu_2)^d$. If the sets $\{u_1=0\}$ and $\{u_2=0\}$ differ from $S$ by a pluripolar set and if $0\leq u_1\leq u_2$ on $\Omega$, then 
$$(dd^cu_1)^d\leq (dd^cu_2)^d.$$
\end{lemma} 

We work in $\CC^2$; indeed, mostly in $\RR^2$. For all of the compact sets $K\subset \RR^2$ in the rest of this paper, $V_K=V_K^*$; i.e., $V_K$ is continuous. The letter $\mathcal C$ denotes a generic constant which can vary from line to line. Let $x_j=\hbox{Re}z_j, \ j=1,2$ where $(z_1,z_2)$ are coordinates for $\CC^2$. We recall Example 4.4 of \cite{SS} which is the union of two quarter-disks meeting at the origin: 
\begin{equation}\label{sset} S:=\{(x_1,x_2)\in \RR^2: x_1^2+x_2^2\leq 1, \ x_1x_2 \geq 0\}.\end{equation}
Let $h(\zeta):=\zeta +\sqrt{\zeta^2-1}$ be the standard Joukowski map in $\CC\setminus [-1,1]$. Sigurdsson and Snaebjarnarson show that
\begin{equation}\label{ssfcn} V_S(z_1,z_2)=\frac{1}{2}\log h(|1-z_1^2-z_2^2|+|z_1-z_2|^2+2|z_1z_2|)\end{equation}
(the formula in \cite{SS} is missing the $1/2$ in front).

An elementary way to obtain (\ref{ssfcn}) is to use Theorem \ref{klimek}. To be precise, the rotation counterclockwise by $45^o$ given by 
$$T_1(x_1,x_2) =\frac{1}{\sqrt 2}(x_1-x_2,x_1+x_2)=(x,y)$$ maps $S$ to $T_1(S)=\tilde S$. Then the square map $Q(x,y)=(x^2,y^2)=(s,t)$ maps $\tilde S$ to 
$$Q(\tilde S)=T=co\{(0,0), (0,1),(1/2,1/2)\}$$
(in a $4$-to-$1$ manner). Here $co(\mathcal S)$ denotes the convex hull of the set $\mathcal S$. Finally, the linear map $T_2(s,t)=(2s,t-s)=(u,v)$ maps $T$ to the standard triangle 
$$\Sigma:=T_2(T)= co\{(0,0), (0,1),(1,0)\}.$$
Altogether, the composition gives the map
$$G(x_1,x_2):=\bigl( (x_1-x_2)^2,2x_1x_2\bigr)=(u,v)$$
with $G^{-1}(\Sigma)=S$. Extending to $\CC^2$, from Theorem \ref{klimek}, 
$$V_{\Sigma}(G(z_1,z_2))=2V_S(z_1,z_2).$$
Using the known formula for $V_{\Sigma}$ (cf., Example 5.4.7 of \cite{K} or Example 4.8 of \cite{Ba}), the left-hand-side equals
$$V_{\Sigma}(G(z_1,z_2))=\log h(|z_1-z_2|^2+2|z_1z_2|+|(z_1-z_2)^2+2z_1z_2-1|)$$
$$=\log h(|z_1-z_2|^2+2|z_1z_2|+|z_1^2+z_2^2-1|).$$
This gives the (correct) formula (\ref{ssfcn}) for $V_S$.

\begin{multicols}{3}

\includegraphics[height=3cm]{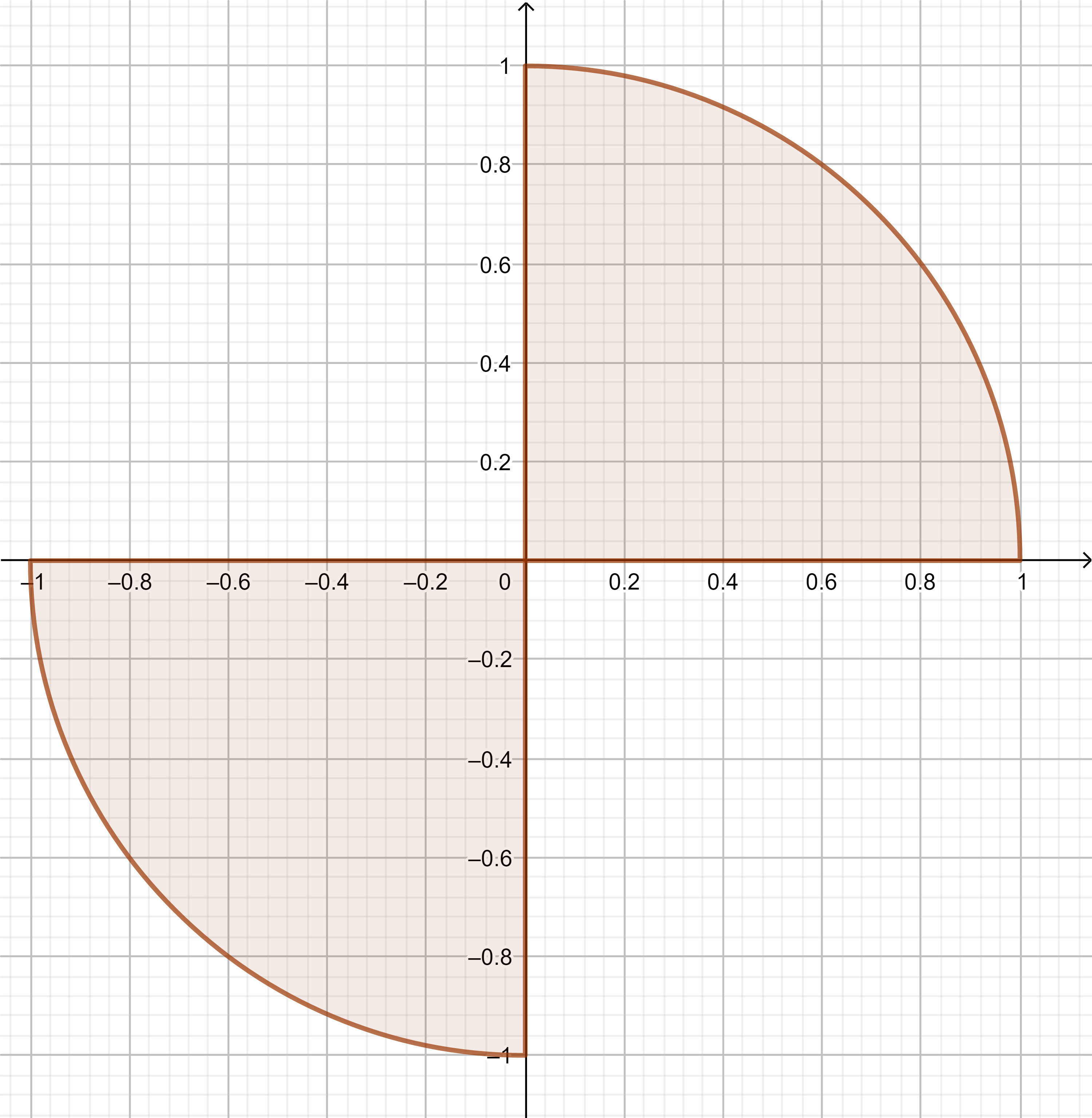} \\
$S$

\includegraphics[height=3cm]{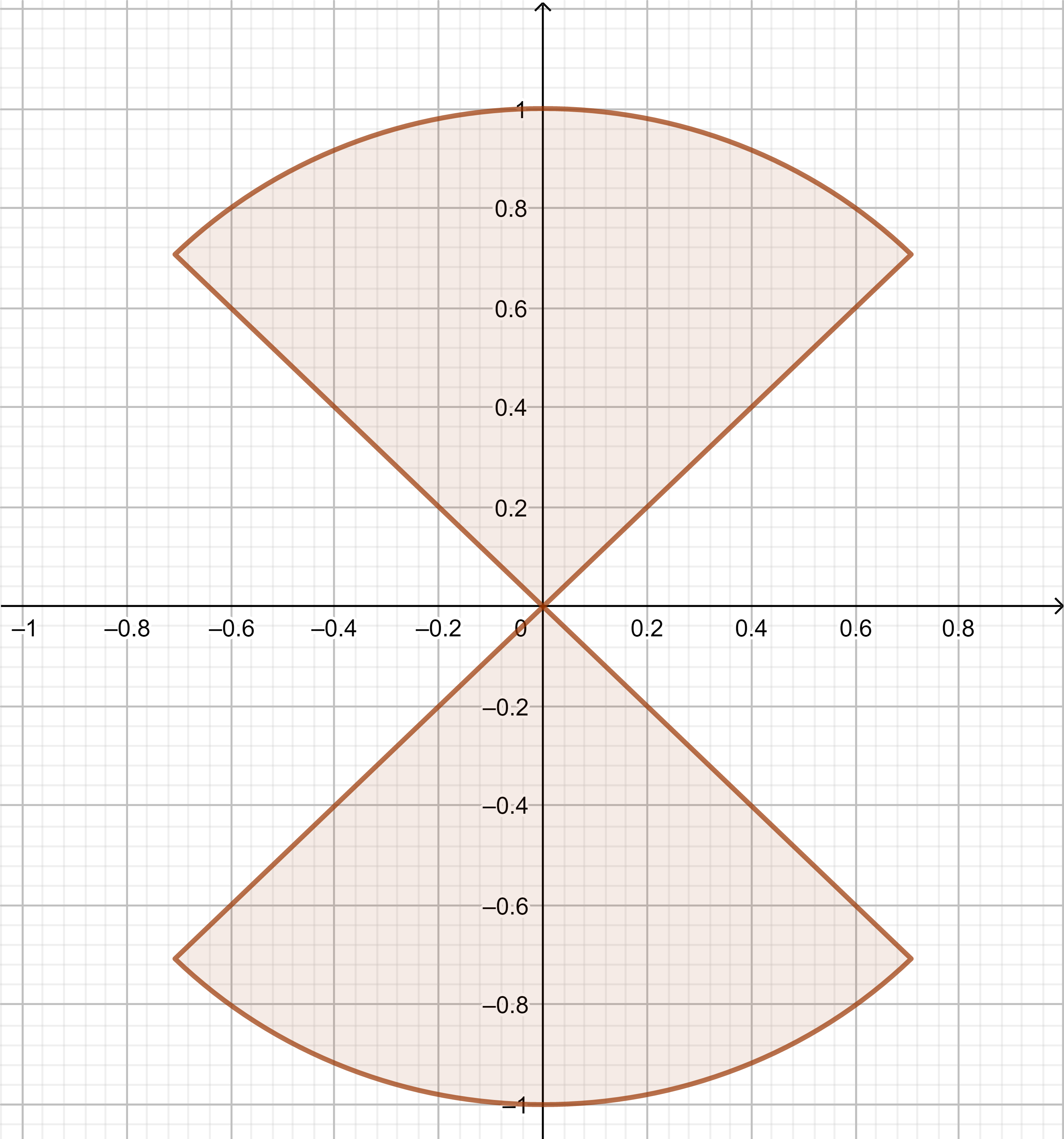} \\
$\tilde S$

\includegraphics[height=3cm]{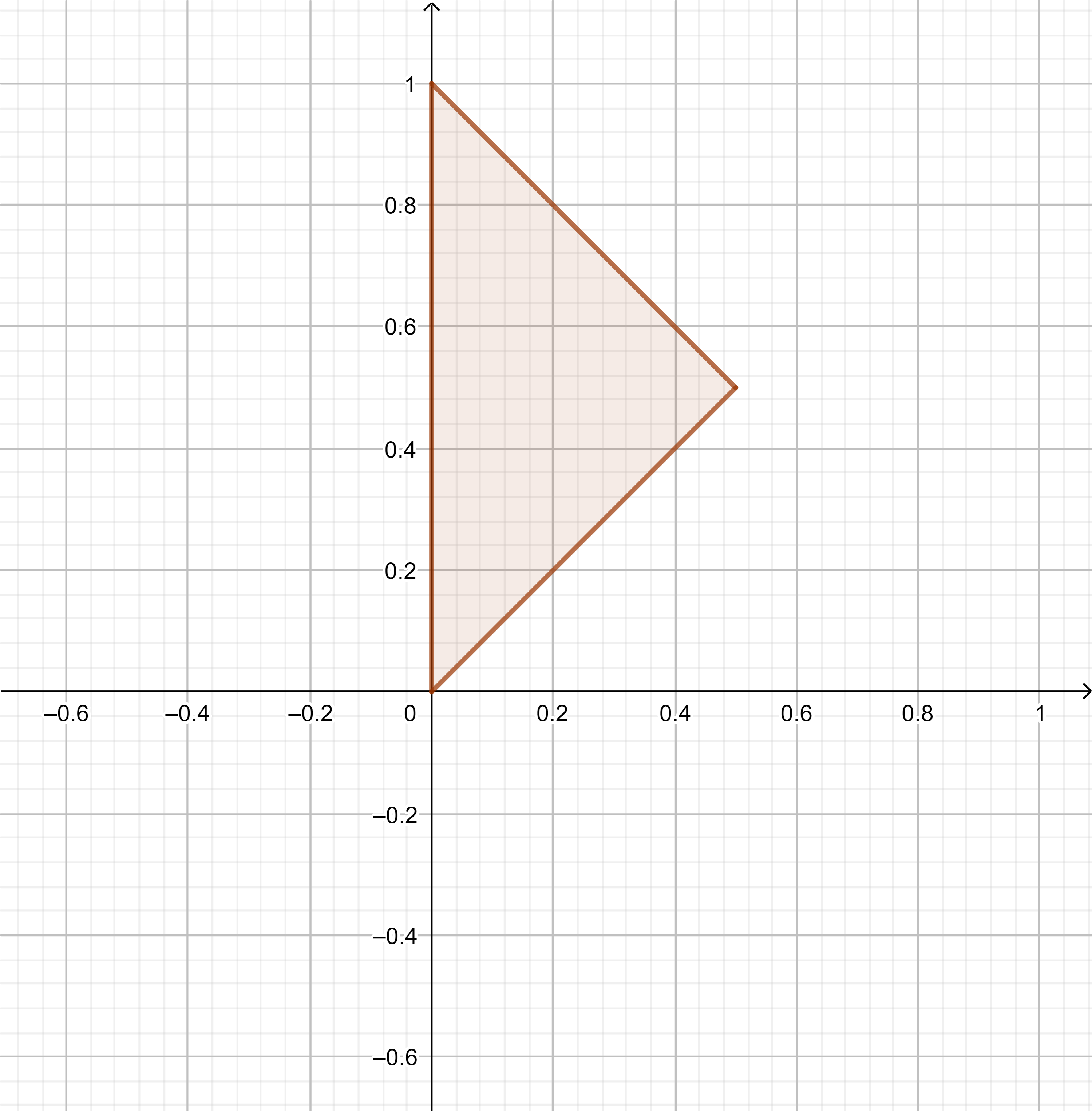} \\
$T$ 
\end{multicols}

More to the point, Sigurdsson and Snaebjarnarson show that the Monge-Amp\`ere density is of the form 
\begin{equation}\label{ssmonge} \rho_S(x_1,x_2)= \mathcal C\frac{|x_1+x_2|}{\sqrt {x_1x_2(1-x_1^2-x_2^2)}}\end{equation}
on the $\RR^2-$interior of $S$ for an appropriate constant $\mathcal C$. Thus as $(x_1,x_2)\to (0,0)$,
$$\rho_S(x_1,x_2)\asymp \frac{|x_1+x_2|}{\sqrt {x_1x_2}}.$$
In particular, along a linear approach $x_2=mx_1$ which is neither vertical nor horizontal, 
$$\lim_{(x_1,x_2)\to (0,0)}\rho_K(x_1,x_2)= \mathcal C \frac{1+m}{\sqrt m}$$
whereas $\rho_S(x_1,x_2)$ blows up as $(x_1,x_2)\to (0,0)$ along a tangential (vertical or horizontal) approach. Note that the linear approach ratio $\frac{1+m}{\sqrt m}$ approaches $+\infty$ as $m\to 0$ or $m\to \infty$; i.e., as we approach an ``edge'' of the quarter disk. In the next section, we show that a similar result holds for an appropriate ``Pac-Man'' set in $\RR^2$.

\section{A Pac-Man set}

In this section we prefer to use $(z,w)$ as coordinates in $\CC^2$ and we let $x=\hbox{Re}z, \ y=\hbox{Re}w$ be the standard coordinates on $\RR^2\subset \CC^2$. Let 
\begin{equation}\label{pacmanp} P=P_{\pi/2}:=\{(x,y)\in \RR^2: (x-1)^2+y^2\leq 1\}\setminus W\end{equation} where 
$$W=W_{\pi/2}:=\{(x,y): x>1, \ |y|<x-1\}$$ is a symmetric wedge of angle opening $\pi/2$ at $(1,0)$; hence the linear boundary of $P$ is contained in the lines $y=\pm(x-1)$ with $x\geq 1$. Thus $P$ is a ``Pac-Man'' of radius one whose ``mouth'' is the wedge $W$. Similar calculations will work for 
$$P_{\alpha}=\{(x,y)\in \RR^2: (x-1)^2+y^2\leq 1\}\setminus W_{\alpha}$$
where $W_{\alpha}$ is a symmetric wedge of angle opening $\alpha <\pi$ at $(1,0)$; for simplicity we take $\alpha=\pi/2$.

Let $F:\CC^2\to \CC^2$ via $F(z,w)=(z,w^2)$. Then $P=F^{-1}(K)$ where $K$ is the {\it convex} set of points in the quadrant $\{(s,t)\in \RR^2: s,t\geq 0\}$ bounded by $\{(s,t): 0\leq s\leq 1, \ t=0\}$, $\{(s,t): t=(s-1)^2, \ 1\leq s \leq 1+1/\sqrt 2\}$ and $\{(s,t): t=2s -s^2, \ 0\leq s \leq 1+1/\sqrt 2\}$. 

From Theorem \ref{klimek}, 
$$V_P(z,w)\leq V_K(F(z,w))\leq 2 V_P(z,w).$$
Thus using Lemma \ref{twoone}, first with $u_1=V_P$ and $u_2 = V_K\circ F$, then with $u_1=V_K\circ F$ and $u_2= 2V_P$, 
\begin{equation}\label{comp}(dd^c V_P)^2\leq \bigl(dd^c (V_K\circ F)\bigr)^2\leq 4 (dd^c V_P)^2.\end{equation}

Note that at $(1,0)\in K$, the tangent line to $\partial K$ is horizontal. Following the argument in the proof of Theorem 1.3 of \cite{BT}, the Monge-Amp\`ere density $\rho_K(s,t)$ of $V_K$ near  $(1,0)$ grows like 
\begin{equation}\label{estimate} \rho_K(s,t)\asymp \frac{1}{dist \bigl((s,t),\partial K\bigr)^{1/2}};\end{equation}
i.e., just as for the Monge-Amp\`ere density of the extremal function for a real disk (or convex body with smooth boundary) or square. Indeed, for $E:=\{(x,y)\in \RR^2: x^2+y^2\leq 1\}$ we have 
$$\rho_E(x,y)=\mathcal C(1-x^2-y^2)^{-1/2}$$ while for 
$\mathcal Q=\{(x,y)\in \RR^2: |x|, |y|\leq 1\}$ we have 
$$\rho_{\mathcal Q}(x,y)=\mathcal C(1-x^2)^{-1/2}(1-y^2)^{-1/2}$$
(cf. \cite{Ba}, Examples 4.6 and 4.7). Moreover, (\ref{estimate}) is a local estimate: to be precise, we first take a neighborhood $N$ of $(1,0)$ constructed as a union of segments $s_{\alpha}$ perpendicular to the tangent line $l_{\alpha}$ to $\partial K$ at $p_{\alpha}$ where the $p_{\alpha}$ fill out a neighborhood of $(1,0)$ in $\partial K$ and $s_{\alpha}$ is sufficiently small but with lengths $|s_{\alpha}|\geq \mathcal C>0$ for some $\mathcal C$ so that we can fit inscribed disks $E_{\alpha}\subset K$ with $p_{\alpha}\in \partial K \cap \partial E_{\alpha}$ of radii $r_{\alpha}\geq |s_{\alpha}|$. Then for points $(s,t)$ on $s_{\alpha}$ we have 
$$dist \bigl((s,t),\partial K \bigr)=dist \bigl((s,t),\partial E_{\alpha} \bigr).$$
Since $V_{E_{\alpha}}\geq V_K$ and the Monge-Amp\`ere density $\rho_{E_{\alpha}}(s,t)$ of $V_{E_{\alpha}}$ behaves like 
$$\rho_{E_{\alpha}}(s,t)\asymp \frac{1}{dist \bigl((s,t),\partial E_{\alpha} \bigr)^{1/2}},$$
we get an upper bound
$$\rho_K(s,t) \leq \frac{\mathcal C}{dist \bigl((s,t),\partial K\bigr)^{1/2}}$$
for $(s,t)\in N$ using Lemma \ref{twoone}. For a lower bound
$$\rho_K(s,t) \geq \frac{\mathcal C}{dist \bigl((s,t),\partial K\bigr)^{1/2}}$$
on $N$ we again use Lemma \ref{twoone} (i.e., the argument in the proof of Theorem 1.3 of \cite{BT}) with circumscribed disks or squares of uniform size (radii or side length bounded {\it above}). To be precise, we let $T_{\alpha}$ be a circumscribed real disk or square (i.e., containing $K$) with $p_{\alpha}\in \partial T_{\alpha}$ and $l_{\alpha}$ tangent to $T_{\alpha}$ at $p_{\alpha}$. Then for points $(s,t)$ on $s_{\alpha}$ 
$$dist \bigl((s,t),\partial T_{\alpha} \bigr)= dist \bigl((s,t),\partial K \bigr).$$
Since $V_{T_{\alpha}}\leq V_K$ and the Monge-Amp\`ere density $\rho_{T_{\alpha}}(s,t)$ of $V_{T_{\alpha}}$ behaves like 
$$\rho_{T_{\alpha}}(s,t)\asymp \frac{1}{dist \bigl((s,t),\partial T_{\alpha} \bigr)^{1/2}},$$
we get the desired lower bound.

We use all these ingredients to show:

\begin{theorem} \label{pac} Let $\rho_P(x,y)$ be the Monge-Amp\`ere density of $V_P$ for $P$ in (\ref{pacmanp}); i.e., 
$$(dd^cV_P)^2=\rho_P(x,y)dx\wedge dy \ \hbox{for} \ (x,y)\in P^o.$$
\begin{enumerate}
\item For an approach to $p=(1,0)$ linearly along $y=c(x-1)$ with $c>1$, 
\begin{equation}\label{limc} \lim_{x\to 1^+}\rho_P(x,c(x-1))= A\frac{|c|}{\sqrt {c^2-1}} \end{equation}
where $A$ is independent of $c$; and for an approach to $p=(1,0)$ linearly along $y=c(x-1)$ with $c\leq 0$ or for a vertical approach along $x=1$ (corresponding to $|c|\to \infty$) we get 
$$\lim_{y\to 0^+}\rho_P(1,y)=A.$$
\item For a tangential approach to $p=(1,0)$; i.e., along $y=(x-1)+g(x)$ with $g(x)\geq 0$ for $x\geq 1$ and $g'(1)=0$, 
\begin{equation}\label{limt} \lim_{x\to 1^+}\rho_P(x,(x-1)+g(x))= \infty. \end{equation}
\end{enumerate}
\end{theorem}

\begin{remark} Note that the linear approach ratio $\frac{|c|}{\sqrt {c^2-1}}$ near the edge $\{(x,y): y=x-1, \ 1\leq x\leq 1+1/\sqrt 2\}$ of the Pac-Man approaches $\infty$ as $c\to 1^+$. For linear approaches ``away'' from the edge, we get the same linear approach value.
\end{remark}

\begin{proof} We begin by observing from (\ref{comp}) that if we write the Monge-Amp\`ere density of $V_P$ as $\rho_P(x,y)$ and that of $V_K$ as $\rho_K(s,t)$ where $(x,y)=(s,t^2)$, 
$$\rho_P(x,y) dx\wedge dy \approx \rho_K(s,t) ds \wedge dt \approx  \rho_K(s,t) \sqrt t dx \wedge dy.$$
Thus, along a curve $y=f(x)$ approaching $p$ as $x\to 1^+$, 
$$\lim_{x\to 1^+}\rho_P(x,f(x))=\lim_{(s,t)\to (1,0), \ \sqrt t =f(s)} [\rho_K(s,t) \cdot \sqrt t ].$$

To prove (\ref{limc}), we consider for $c>1$
\begin{equation}\label{help}\lim_{x\to 1^+}\rho_P(x,c(x-1))=\lim_{(s,t)\to (1,0), \ \sqrt t =c(s-1)} [\rho_K(s,t) \cdot \sqrt t ].\end{equation}
Now along $\sqrt t =c(s-1)$ as $(s,t)\to (1,0)$, we have 
$$dist \bigl((s,t),\partial K\bigr)\asymp c^2(s-1)^2-(s-1)^2=(c^2-1)(s-1)^2$$
(this is the vertical distance between points $\sqrt t =c(s-1)$ and $\sqrt t =(s-1)$). Hence together with (\ref{estimate}) we obtain
$$\rho_K(s,t) \asymp \frac{1}{dist \bigl((s,t),\partial K \bigr)^{1/2}} \asymp \frac{1}{(s-1)\sqrt {c^2-1}}= \frac{c}{\sqrt t \sqrt {c^2-1}}.$$
Thus, from (\ref{help}),
$$\lim_{x\to 1^+}\rho_P(x,c(x-1))= A\frac{c}{\sqrt {c^2-1}}$$
for some constant $A$ which is independent of path. For $c<0$ 
$$dist \bigl((s,t),\partial K\bigr) =t $$
since the closest point to $(s,t)$ in $\partial K$ is $(s,0)$. Thus
$$\rho_K(s,t) \asymp \frac{1}{dist \bigl((s,t),\partial K \bigr)^{1/2}} \asymp \frac{1}{\sqrt t}$$
and from (\ref{help}), which is still valid for such $c$, 
$$\lim_{x\to 1^+}\rho_P(x,c(x-1))= A.$$ 
Similarly, along $x=1$, we have $s=1$ and the (vertical) distance between points $(1,t)$ and $\partial K$ (with closest point $(1,0)$) is $ t$ so $\rho_K(s,t) \cdot \sqrt t \asymp 1/\sqrt t \cdot \sqrt t =1$ and thus
$$\lim_{(1,t)\to (1,0)} [\rho_K(s,t) \cdot \sqrt t ]=A.$$

To prove (\ref{limt}), we consider 
$$\lim_{x\to 1^+}\rho(x,(x-1)+g(x))=\lim_{(s,t)\to (1,0), \ \sqrt t =(s-1)+g(s)} [\rho_K(s,t) \cdot \sqrt t ].$$
Note that both the horizontal distance $H(s,t)$ from a point $(s,t)$ with $\sqrt t =(s-1)+g(s)$ to $\partial K$ and the vertical distance $V(s,t)$ are greater than or equal to  $dist \bigl((s,t),\partial K\bigr)$. Now $H(s,t)=g(s)$; if instead we consider the vertical distance $V(s,t)$ we have the estimate
$$V(s,t)= [(s-1)+g(s)]^2-(s-1)^2=2(s-1)g(s) +g(s)^2\geq dist \bigl((s,t),\partial K\bigr).$$
Thus
$$\rho_K(s,t) \cdot \sqrt t \asymp \frac{ \sqrt t }{dist \bigl((s,t),\partial K \bigr)^{1/2}} \geq \frac{ \sqrt t }{V(s,t)^{1/2}}$$
$$=\bigl(\frac{[(s-1)+g(s)]^2}{[2(s-1)g(s) +g(s)^2]}\bigr)^{1/2}=\bigl(1+\frac{(s-1)^2}{[2(s-1)g(s) +g(s)^2]}\bigr)^{1/2}.$$
If we expand $g(s)=a(s-1)^N+ ...$ where $N> 1$ we have 
$$\rho_K(s,t) \cdot \sqrt t \geq 0\bigl((1+  \frac{(s-1)^2}{2a(s-1)^{N+1}})^{1/2}\bigr)\to \infty$$
as $s\to 1^+$ since $N>1$. This proves (\ref{limt}).
\end{proof}

\begin{remark} As previously mentioned, similar calculations work for other Pac-Men $P_{\alpha}$ as long as the opening angle $\alpha$ is less than $\pi$. If the angle is greater than or equal to $\pi$, $P_{\alpha}$ is {\it convex} and the Monge-Amp\`ere density of $V_{P_{\alpha}}$ becomes unbounded near $\partial P_{\alpha}$ regardless of the direction of approach. 

\end{remark} 

\begin{remark} The same behavior holds for a fattened $L-$shaped region at its inner bend $p$ as can be seen by putting an inscribed Pac-Man with vertex $p$ inside $L$ and a circumscribed Pac-Man with vertex $p$ containing $L$ and using Lemma \ref{twoone}. For example, setting $L=E\cup E^*$ where 
$$E:=\{(x,y)\in \RR^2: 0\leq y\leq 1+\sqrt 2 -x, \ x-1\leq y \leq x+\sqrt 2 - 1\}$$
and $E^*=\{(x,-y):(x,y)\in E\}$, the Pac-Man $P$ in (\ref{pacmanp}) is inscribed in $L$ with $p=(1,0)$.

\end{remark}



\begin{remark} \label{ragnar} We recover a (rotated, translated) version $S'$ of the two quarter disk set $S$ in (\ref{sset}) by taking out a symmetric wedge $W'$ from $P$; thus the linear boundary of $S'$ is contained in the lines $y=\pm(x-1)$ with $1-1/\sqrt 2 \leq x \leq 1+1/\sqrt 2$. Here $S'=F^{-1}(K')$ where, as before, $F(z,w)=(z,w^2)$, and $K'$ is now the convex set bounded by the parabolas $\{(s,t): t=(s-1)^2\}$ and $\{(s,t): t=2s -s^2\}$. Thus one can follow the proof of Theorem \ref{pac} to verify that we do, indeed, have the same asymptotic results for the density $\rho_{S'}$ for $S'$ near $(1,0)$ and for the density $\rho_S$ for
$$S=\{(u,v)\in \RR^2: 0\leq u^2+v^2\leq 1, \ uv\geq 0\}$$ 
near $(0,0)$. (Here $S$ is the same set as in (\ref{sset}) but we are using coordinates $(u,v)=(x_1,x_2)$.) From (\ref{ssmonge}) the Monge-Amp\`ere density of $V_S$ is $\rho_S(u,v)\asymp \frac{|u+v|}{\sqrt {uv}}$; hence, e.g., along lines $v=mu$ for $0< m < 1$ and $u,v\geq 0$, we obtain
$$\lim_{(u,v)\to (0^+,0^+), \ v=mu}\rho_S(u,v)= B(\frac{1+m}{m})$$
for some constant $B$ independent of $m$. One easily checks, via translation and rotation of coordinates, that this corresponds to $$1< c=\frac{m+1}{1-m} <\infty \ \hbox{or} \ \frac{c}{\sqrt {c^2-1}}=\frac{1}{4}\cdot \frac{1+m}{m}$$ where $c$ is as in Theorem \ref{pac}.

\end{remark}

\section{Pluripotential theory and convex bodies}

The impetus for our proof of Theorem \ref{pac} came from the newly-developed pluripotential theory associated to convex bodies (cf., \cite{BBL}). For $C\subset (\RR^+)^d$ a convex body define  
$$H_C(z):=\sup_{J\in C} \log |z^J|:=\sup_{(j_1,...,j_d)\in P} \log[|z_1|^{j_1}\cdots |z_d|^{j_d}],$$
the logarithmic indicator function of $C$. Thus for $\Sigma:=\{(x_1,...,x_d)\in \RR^d: 0\leq x_i \leq 1, \ \sum_{j=1}^d x_i \leq 1\}$, 
$H_{\Sigma}(z)=\max_{j=1,...,d}\log^+ |z_j|$. We assume $\Sigma \subset kC \ \hbox{for some} \ k\in \ZZ^+  $. Define
$$L_C=L_C(\CC^d):= \{u\in PSH(\CC^d): u(z)- H_C(z) =0(1), \ |z| \to \infty \} \ \hbox{and}$$ 
$$Poly(nC):=\{p(z)=\sum_{J\in nC\cap (\ZZ^+)^d}c_J z^J: c_J \in \CC\}, \ n=1,2,...$$
For $p\in Poly(nC)$ we have $\frac{1}{n}\log |p|\in L_C$. Note $Poly(n\Sigma)$ are the usual holomorphic polynomials of degree at most $n$.

Given $E\subset \CC^d$, the {\it $C-$extremal function of $E$} is given by $V^*_{C,E}(z)$ where
$$V_{C,E}(z):=\sup \{u(z):u\in L_C(\CC^d), \ u\leq 0 \ \hbox{on} \ E\}.$$
For $K\subset \CC^d$ compact,
$$V_{C,K}(z)=\lim_{n\to \infty} [\sup  \{  \frac{1}{n}\log |p(z)|: p\in Poly(nC), \ ||p||_K\leq 1\}].$$
Thus $V_{\Sigma,K}=V_K$ from the introduction. Moreover, following the proof of Theorem 5.3.1 of \cite{K}, we have the following result.

\begin{theorem}\label{cklimek} Let $C,C'\subset(\RR^+)^d$ be convex bodies and let $F:\CC^d\to \CC^d$ be a proper polynomial mapping satisfying
$$0<\liminf_{|z|\to \infty}  \frac{\sup_{J\in C} |[F(z)]^{J}|}{\sup_{J\in C'} |z^{J'}|} \leq \limsup_{|z|\to \infty}  \frac{\sup_{J\in C} |[F(z)]^J|}{\sup_{J'\in C'} |z^{J'}|} < \infty.$$
Then
$$V_{C,K}(F(z))=V_{C',F^{-1}(K)}(z).$$

\end{theorem} 

Relating to the proof of Theorem \ref{pac}, we observe that for the mapping $F:\CC^2\to \CC^2$ given by $F(z,w)=(z,w^2)$ that we used, given {\it any} $P'=F^{-1}(K')\subset \{(x,y):(x-1)^2+y^2\leq 1\}$ and the corresponding set $K'$, we have 
$$2V_{P'}(z,w)= V_{2\Sigma,P'}(z,w)=V_{C,K'}(F(z,w))$$
where $C$ is the convex hull of $(0,0), (2,0)$ and $(0,1)$. This follows since $F$ satisfies the hypothesis of Theorem \ref{cklimek}. 

We hope that one can use this elementary result to construct other explicit examples of these $C-$extremal functions, $V_{C,K}$. The most useful case(s) would be when one of $C$ or $C'$ is $\Sigma$.


\bigskip

\begin{thebibliography}{GGK2}
\bibitem{Ba} M. Baran, Complex equilibrium measure and Bernstein type theorems for compact sets in $\RR^n$, \emph{Proc. AMS}, \textbf{123}, (1995), no. 2, 485-494.
\bibitem{BT} E. Bedford and B. A. Taylor, The complex equilibrium measure of a symmetric convex set in $\RR^n$, \emph{Trans. AMS}, \textbf {294}, (1986), no. 2, 705-717.
\bibitem{BBL} T. Bayraktar, T. Bloom and N. Levenberg, Pluripotential theory and convex bodies, \emph{Mat. Sbornik}, \textbf{209}, (2018), no. 3, 352-384.


\bibitem{K} M. Klimek, \emph{Pluripotential Theory}, Oxford University Press, 1991.

\bibitem{SS} Ragnar Sigurdsson, Audunn Skuta Snaebjarnarson, Monge-Amp\`ere measures of plurisubharmonic exhaustions associated to the Lie norm of holomorphic maps, arXiv:1810.01326 
 

\end{thebibliography}
\end{document}